\documentclass{amsart}
\usepackage{amssymb}
\usepackage{mathrsfs}
\usepackage{graphicx}
\usepackage{comment}
\newcommand{\comments}[1]{}
\usepackage[francais,english]{babel}

\def\midd{\, : \,}
\newcommand{\om}{\omega}

\newcommand{\be}{\begin{equation}}
\newcommand{\ee}{\end{equation}}
\newcommand{\ba}{\begin{align}}
\newcommand{\ea}{\end{align}}

\newtheorem{theorem}{Theorem}[section]

\def\SS{\mathcal S}
{\begin{list}{}{%
\settowidth{\labelwidth}{\textsf{{\it #1.}}}%
\setlength{\labelsep}{4mm}%
\setlength{\leftmargin}{\labelwidth}%
\addtolength{\leftmargin}{\labelsep}%
}}%
{\end{list}}

{\begin{list}{}{%
\settowidth{\labelwidth}{\textsf{{\it #1.}}}%
\setlength{\labelsep}{2mm}%
\setlength{\leftmargin}{\labelwidth}%
\addtolength{\leftmargin}{\labelsep}%
\addtolength{\leftmargin}{4mm}%
\setlength{\itemsep}{6pt}%
\setlength{\listparindent}{0pt}%
\setlength{\topsep}{3pt}%
}}%
{\end{list}}

\title[Integer group determinants for  $S_4$]{The  integer group determinants for the symmetric group of degree four}

\author[C. Pinner]{Christopher Pinner}
\address{ Department of Mathematics\\
         Kansas State University\\
         Manhattan, KS 66506, USA}
\email{pinner@math.ksu.edu}
\thanks{The author thanks Chris Smyth for issuing the challenge of $S_4$.}

\date{\today}
\begin{document}

\begin{abstract}
For the symmetric group $S_4$ we determine all the integer values taken  by its
group determinant when the matrix entries are integers.

\end{abstract}

\maketitle

\section{Introduction}\label{intro}

Suppose that $G=\{g_1,\ldots ,g_n\}$ is a finite group.  Assigning variables $x_g$, $g\in G$, one defines the
{\em group determinant,} $\mathscr{D}_G(x_{g_1},\ldots ,x_{g_{n}})$, to be  the determinant of the $n\times n$ group matrix whose $ij$th entry is  $x_{g_i g_j^{-1}}$. Plainly $\mathscr{D}_G(x_{g_1},\ldots ,x_{g_{n}})$ will be a homogeneous polynomial of degree $n$ in the $x_g$. 

 For $G=\mathbb Z_n$, the cyclic group of order $n$,  this is a circulant determinant, where the next row in the group matrix is obtained from the previous row by a cyclic shift of one to the right. 
An old problem of Olga Taussky-Todd is to determine which  integers can be achieved as an $n\times n$ circulant determinant when the entries in the matrix are integers. We can of course ask the same question for an arbitrary finite  group $G$ and 
define $\SS(G)$ to be this set of integers:
$$ \SS (G) :=\{ \mathscr{D}_G(x_{g_1},\ldots ,x_{g_n})  \midd x_{g_1},\ldots ,x_{g_n}\in \mathbb Z\}. $$
Notice that $\SS(G)$ is closed under multiplication
$$ \mathscr{D}_G(a_{g_1},\ldots ,a_{g_n})\mathscr{D}_G(b_{g_1},\ldots ,b_{g_n})=\mathscr{D}_G(c_{g_1},\ldots ,c_{g_n}), \;\;  \; c_{g}:=\sum_{\substack{u,v\in G\\ uv=g}} a_ub_v,$$
corresponding to multiplication $\left(\sum_{g\in G} a_g g\right)\left(\sum_{g\in G} b_g g\right) =\sum_{g\in G} c_g g$ in the group ring.
While it has been proved \cite{Formanek} that the group determinant polynomial determines the group it is not known whether the set of integer values $\SS (G)$ determines the group.

For $G=\mathbb Z_n$  Newman \cite{Newman1} and Laquer\cite{Laquer} gave divisibility restrictions on the elements in $\SS (G)$ and some  attainable values. For example, $n^2$ and any $m$ with $\gcd(m,n)=1$ 
will be in $\SS(\mathbb Z_n)$, and if $m$ is in $\SS(\mathbb Z_n)$ 
then so is $-m$ and if $p\mid m$ and $p^{\alpha}\parallel n$ then $p^{\alpha+1}\mid m$. Here and throughout $p$ denotes a prime. Newman and Laquer both  obtained a complete description for  $G=\mathbb Z_p$  when $p\geq 2$, and Laquer for $G=\mathbb Z_{2p}$ when $p\geq 3$:
\begin{align*}
\SS (\mathbb Z_p) & = \{ p^a m\midd \gcd(m,p)=1,\;\; a=0 \hbox{ or } a\geq 2\},\\
\SS (\mathbb Z_{2p}) & = \{ 2^a p^b m\midd \gcd(m,2p)=1,\;\; a =0 \hbox{ or } a\geq 2,\;\; b =0 \hbox{ or } b\geq 2\},
\end{align*}
cases where the divisibility conditions are necessary and sufficient.
Newman \cite{Newman2} similarly showed that 
$$\SS(\mathbb Z_{9})=\{3^am\midd \gcd(m,3)=1,\; a=0 \hbox{ or } a\geq 3\},$$
but that $p^3$ is not in $\SS(\mathbb Z_{p^2})$ for any  $p\geq 5$ (in particular the basic divisibility 
conditions are not sufficient). Currently no $\SS(\mathbb Z_{p^2})$ with $p\geq 5$  has  been fully determined, though there are upper and lower set inclusions.

As shown by Vipismakul \cite{Cid1} there is a very close  relationship between the group determinant for a finite abelian group and Lind's generalization  \cite{Lind} of the Mahler measure  for that group; in particular the corrresponding  Lind-Lehmer problem of  finding the minimal positive logarithmic 
Lind-Mahler  measure for the group,   corresponds to finding the smallest non-trivial group determinant
$$ \lambda(G) :=\min \{ |m| \midd m\in \SS (G), \; |m|\geq 2\}. $$
Kaiblinger \cite{Norbert2} used the Lind measure approach to find 
$$\SS(\mathbb Z_4)=\{ 2m+1,\;\; 2^4m\midd m\in \mathbb Z\}, $$ 
and $\SS(\mathbb Z_8),$ and obtain
upper and lower inclusions for other $\SS(\mathbb Z_{2^k})$. The Lind-Lehmer constant  $\lambda (G)$ is known for a number of groups, see for example   \cite{Norbert,Pigno1,pgroups,Cid2}, including the cyclic groups $G=\mathbb Z_n$ with $892371480\nmid n$.
Writing $D_{2n}$ for the dihedral group of order $2n$,  the sets 
$\SS(D_{2p})$, $\SS(D_{4p}),$ 
$$\SS(D_{8})=\{4m+1,\;\; 2^8m \midd m\in \mathbb Z\}, $$
and $\SS(D_{16})$ were determined in \cite{dihedral}, and $\lambda(D_{2n})$ found for $n< 1.89\times 10^{47}$.
In particular the cases $p=2$ and $3$ give us  the Klein $4$-group
$$ \SS(\mathbb Z_2\times \mathbb Z_2)=\{4m+1,\; \; 2^4(2m+1),\; \; 2^6m\midd m\in \mathbb Z\}, $$
and the symmetric group of degree $3$
$$   \SS(S_3)= \{  2^a3^bm \midd \gcd(m,6)=1,\; a=0 \hbox{ or } a\geq 2,\; b=0 \hbox{ or } b\geq 3    \}. $$

In \cite{smallgps} a complete description of the integer group determinants was obtained for the remaining groups of order at most $14$, including   the alternating group $A_4$:
$$\SS(A_4)=\SS(A_4)_{odd} \cup \SS(A_4)_{even},$$ with
\begin{align*} \SS(A_4)_{odd} & =\{m\equiv 1 \text{ mod } 4 \midd 3\nmid m \hbox{ or } 3^2\mid m\},\\
\SS(A_4)_{even} & =  \{ 2^a3^bm\midd \gcd(m,6)=1,\; a=4 \hbox{ or } a\geq 8, \; b=0 \hbox{ or } b\geq 2\}. \end{align*}
Although $\SS(S_3)$ and $\SS(A_4)$ were obtained without too much difficulty, other small groups  in \cite{smallgps}, for example $G= \mathbb Z_2\times \mathbb Z_6$, were considerably more complicated; making it clear  that obtaining a description for  general  $\SS (G)$ is probably not feasible. Indeed, even in the case of circulant determinants,  we are yet to obtain a complete description   of $\SS (\mathbb Z_{15})$ and $\SS(\mathbb Z_{16})$. Our goal here is to show that $G=S_4$ is one of those rare cases where we can completely  determine $\SS(G)$:

\begin{theorem}\label{main}
For $G=S_4$ the odd group  determinants in $\SS(G)$  are the integers $m=1$ mod $4$ with the property that $3\nmid m$ or $3^3\mid m$.

The even determinants are the  integers of the form
$$ 2^8m, \;\; m\equiv 1 \text{ mod }4,\;\;\;  \hbox{ or }\; \;\;2^{10}m, \;\; m\equiv -1 \text{ mod }4,\;\;\; \hbox{ or }\;\;\;  2^{12}m ,\;\;  $$
where $m$ is an integer with $3\nmid m$ or $3^3\mid m$.

\end{theorem}
Comparing  $\SS(S_4)$ to $\SS(G)$ for $G=\mathbb Z_2,\mathbb Z_3, \mathbb Z_4,\mathbb Z_2\times \mathbb Z_2,D_8,S_3,A_4$ 
it is  tempting to ask whether in general $\SS(G)\subseteq \SS(H)$ whenever $H\leq G$.

\section{Factoring the group determinant for  $S_4$} 

Dedekind \cite{Dedekind} observed that for a finite abelian group the group characters $\hat{G}$ couuld be used to factor
the group determinant into linear factors
\be \label{abelian}  \mathscr{D}_G(x_{g_1},\ldots ,x_{g_n}\}=\prod_{\chi \in \hat{G}} \left(\chi(g_1) x_{g_1}+\cdots +\chi(g_n)x_{g_n}\right). \ee
For non-abelian groups the group determinant will contain non-linear factors and, as  shown by Frobenius \cite{Frobenius}, the
counterpart to \eqref{abelian} is to  use  the set $\hat{G}$ of irreducible  group representations for $G$,
$$  \mathscr{D}_G(x_{g_1},\ldots ,x_{g_n}\}=\prod_{\rho \in \hat{G}}\det  \left(x_{g_1}\rho(g_1)+\cdots +x_{g_n}\rho(g_n)\right)^{\deg(\rho)},    $$
see Conrad \cite{Conrad} for a historical survey.  

The irreducible representations for $S_4$ are discussed in Serre \cite[$\S$5.8]{Serre}. One can also use GAP or a similar computer algebra system to generate them  (although we have reduced the four generators used  there to two and also reindexed). 

We take as our two generators for $S_4$
$$ \alpha=(1234),\;\;\;\beta=(12), $$
and order the even permutations, with the coefficients $a_1,\ldots a_{12}$, by
\begin{align*} & \{ 1,(13)(24),(14)(23),(12)(34),(134),(243),(142),(123),(143),(132),(124),(234)\}\\
 & =\{1,\alpha^2,\beta\alpha^2\beta,\alpha^2\beta\alpha^2\beta,\alpha\beta,\alpha^3\beta,\alpha\beta\alpha^2,\alpha^3\beta\alpha^2,\beta\alpha^3,\alpha^2\beta\alpha,\alpha^2\beta\alpha^3,\beta\alpha\}\end{align*}
and the odd permutations, with coefficients $b_1,\ldots ,b_{12},$ by
\begin{align*} & \{ (1234),(1432),(24),(13),(14),(23),(1243),(1342),(12),(34),(1324),(1423)\}\\
 & = \{ \alpha,\alpha^3,\alpha\beta\alpha^2\beta,\beta\alpha^2\beta\alpha,\alpha^3\beta\alpha,\alpha\beta\alpha^3,\alpha\beta\alpha,\alpha^3\beta\alpha^3,\beta,\alpha^2\beta\alpha^2,\beta\alpha^2,\alpha^2\beta\}.
\end{align*}
For $S_4$ we have two linear representations, $\chi_0(x)=1$ and $\chi_1(x)=\hbox{sgn}(x)$, giving two linear factors
\begin{align*}
\ell_1:=(a_1+\cdots + a_{12})+(b_1+\cdots +b_{12}),\\
\ell_2:=(a_1+\cdots + a_{12})-(b_1+\cdots +b_{12}).
\end{align*}
We have one degree two representation with 
$$ \rho_1(\beta)=\left(\begin{matrix} 0 & \om^2 \\ \om & 0\end{matrix}\right),\;\;\;\rho_1(\alpha)=\left(\begin{matrix} 0 & 1 \\ 1 & 0\end{matrix}\right),\;\;\;    \om := e^{2\pi i/3}, $$
giving
\begin{align*}  \rho_1(1)  =\rho_1(\alpha^2)=\rho_1(\beta\alpha^2\beta)=\rho_1(\alpha^2\beta\alpha^2\beta) &=\left(\begin{matrix} 1 & 0 \\ 0 & 1 \end{matrix}\right), \\
\rho_1(\alpha\beta) = \rho_1(\alpha^3\beta)=\rho_1(\alpha\beta\alpha^2)=\rho_1(\alpha^3\beta\alpha^2)&=\left(\begin{matrix} \om & 0 \\ 0 & \om^2 \end{matrix}\right), \\
\rho_1(\beta\alpha^3) =\rho_1(\alpha^2\beta\alpha)=\rho_1(\alpha^2\beta\alpha^3)=\rho_1(\beta\alpha)&=\left(\begin{matrix} \om^2 & 0 \\ 0 & \om \end{matrix}\right), \\
 \rho_1(\alpha)=\rho_1(\alpha^3)=\rho_1(\alpha\beta\alpha^2\beta)=\rho_1(\beta\alpha^2\beta\alpha)&=\left(\begin{matrix} 0 & 1 \\ 1 & 0 \end{matrix}\right), \\
\rho_1(\alpha^3\beta\alpha)=\rho_1(\alpha\beta\alpha^3)=\rho_1(\alpha\beta\alpha)=\rho_1(\alpha^3\beta\alpha^3)&=\left(\begin{matrix} 0 & \om \\ \om^2 & 0 \end{matrix}\right), \\
\rho_1(\beta)=\rho_1(\alpha^2\beta\alpha^2)=\rho_1(\beta\alpha^2)=\rho_1(\alpha^2\beta)&=\left(\begin{matrix} 0 & \om^2 \\ \om & 0 \end{matrix}\right), 
\end{align*}
and, writing
\begin{align*}
u_1 & =a_1+a_2+a_3+a_4,\;\;\; u_2=a_5+a_6+a_7+a_8,\;\;\; u_3=a_9+a_{10}+a_{11}+a_{12},\\
v_1 & =b_1+b_2+b_3+b_4,\;\;\; v_2=b_5+b_6+b_7+b_8,\;\;\; v_3=b_9+b_{10}+b_{11}+b_{12},
\end{align*}
the quadratic factor
\begin{align*}  q_1:& = \det \left(\begin{matrix} u_1+u_2\;\om+u_3\; \om^2 & v_1+v_2\;\om +v_3\;\om^2 \\ v_1+v_2\;\om^2+v_3\;\om & u_1+u_2\;\om^2+u_3\;\om \end{matrix} \right) \\
 & =N(u_1+u_2\; \om +u_3\; \om^2) - N(v_1+v_2\; \om +v_3\; \om^2),
\end{align*}
 where $N(x)=|x|^2$.

Finally we have two degree three representations
$$ \rho_2(\alpha)=  \left( \begin{matrix}   0 & -1 & 0 \\ 1 & 0 & 0 \\ 0 & 0 & -1 \end{matrix} \right),\;\;\;    \rho_2(\beta)= \left( \begin{matrix}1 & 0 & 0 \\ 0 & 0 & 1 \\ 0 & 1 & 0    \end{matrix} \right), $$
and $\rho_3(x)=\hbox{sgn}(x)\rho_2(x)$. Here $\rho_2$ comes from the natural representation of $S_4$ in $\mathbb C^3$. That is, we take the $3$-dimensional subspace 
$$V=\{x_1e_1+x_2e_2+x_3e_3+x_4e_4\midd x_1+x_2+x_3+x_4=0\} $$
of a  $4$-dimensional vector space and let $S_4$ permute its  basis vectors $e_1,e_2,e_3, e_4$, although to produce  our  $\rho_2$  we take a  less obvious  basis for $V$:
$$ e_1+e_2-e_3-e_4,\;\; -e_1+e_2+e_3-e_4,\;\; e_1-e_2+e_3-e_4. $$

 For the even permutations we have
\begin{align*}
  \rho_2(1) &=\left(\begin{matrix} 1 & 0 & 0 \\ 0 & 1 & 0 \\ 0 & 0 & 1 \end{matrix}\right), &
\rho_2((13)(24))&=\left(\begin{matrix} -1 & 0 & 0 \\ 0 & -1 & 0 \\ 0 & 0 & 1 \end{matrix}\right), &
\rho_2((14)(23))&=\left(\begin{matrix} -1 & 0 & 0 \\ 0 & 1 & 0 \\ 0 & 0 & -1 \end{matrix}\right),\\
 \rho_2((12)(34))& =\left(\begin{matrix} 1 & 0 & 0 \\ 0 & -1 & 0 \\ 0 & 0 & -1 \end{matrix}\right), &
\rho_2((134)) &=\left(\begin{matrix} 0 & 0 & -1 \\ 1 & 0  & 0 \\ 0 & -1 & 0 \end{matrix}\right), &
\rho_2((243)) &=\left(\begin{matrix} 0 & 0 & 1 \\ -1 & 0  & 0 \\ 0 & -1 & 0 \end{matrix}\right),\\
\rho_2((142))& =\left(\begin{matrix} 0 & 0 & -1 \\ -1 & 0  & 0 \\ 0 & 1 & 0 \end{matrix}\right),&
\rho_2((123)) &=\left(\begin{matrix} 0 & 0 & 1 \\ 1 & 0  & 0 \\ 0 & 1 & 0 \end{matrix}\right), &
\rho_2((143)) &=\left(\begin{matrix} 0  & 1 & 0 \\ 0  & 0 & -1 \\ -1 & 0 & 0 \end{matrix}\right),\\
 \rho_2((132)) &=\left(\begin{matrix} 0  & 1 & 0 \\ 0  & 0 & 1 \\ 1 & 0 & 0 \end{matrix}\right),&
\rho_2((124)) &=\left(\begin{matrix} 0  & -1 & 0 \\ 0  & 0 & 1 \\ -1 & 0 & 0 \end{matrix}\right),&
\rho_2((234)) &=\left(\begin{matrix} 0  & -1 & 0 \\ 0  & 0 & -1 \\ 1 & 0 & 0 \end{matrix}\right),
\end{align*}
and for the odd permutations
\begin{align*}
\rho_2((1234))&=\left( \begin{matrix}   0 & -1 & 0 \\ 1 & 0 & 0 \\ 0 & 0 & -1 \end{matrix} \right), &
\rho_2((1432))&=\left( \begin{matrix}   0 & 1 & 0 \\ -1 & 0 & 0 \\ 0 & 0 & -1 \end{matrix} \right), &
\rho_2((24))&=\left( \begin{matrix}   0 & -1 & 0 \\ -1 & 0 & 0 \\ 0 & 0 & 1 \end{matrix} \right), \\
\rho_2((13))&=\left( \begin{matrix}   0 & 1 & 0 \\ 1 & 0 & 0 \\ 0 & 0 & 1 \end{matrix} \right), &
\rho_2((14))&=\left( \begin{matrix}   0 & 0 & -1 \\  0 & 1 & 0  \\ -1 & 0 & 0  \end{matrix} \right), &
\rho_2((23))&=\left( \begin{matrix}   0 & 0 & 1 \\  0 & 1 & 0  \\ 1 & 0 & 0  \end{matrix} \right), \\
\rho_2((1243)) &=\left( \begin{matrix}   0 & 0 & 1 \\  0 & -1 & 0  \\ -1 & 0 & 0  \end{matrix} \right),&
\rho_2((1342)) &=\left( \begin{matrix}   0 & 0 & -1 \\  0 & -1 & 0  \\ 1 & 0 & 0  \end{matrix} \right),&
\rho_2((12)) &=\left( \begin{matrix}  1 &  0 & 0  \\  0 & 0 & 1  \\  0 & 1 & 0   \end{matrix} \right),\\
\rho_2((34)) &=\left( \begin{matrix}  1 &  0 & 0  \\  0 & 0 & -1  \\  0 & -1 & 0   \end{matrix} \right), & 
\rho_2((1324)) &=\left( \begin{matrix}  -1 &  0 & 0  \\  0 & 0 & 1  \\  0 & -1 & 0   \end{matrix} \right), &
\rho_2((1423)) &=\left( \begin{matrix}  -1 &  0 & 0  \\  0 & 0 & -1  \\  0 & 1 & 0   \end{matrix} \right). 
\end{align*}

Thus we have two cubic factors
$$ d_1:=\det \left( A+B\right),\;\;\; d_2:=\det\left(A-B\right), $$
where
$$ A=\left(\begin{matrix}    a_1-a_2-a_3 +a_4 & a_9+a_{10}-a_{11}-a_{12}  & -a_5   +a_6-a_7+a_8    \\  a_5-a_6-a_7+a_8     &  a_1-a_2+a_3-a_4 & -a_9+a_{10}+a_{11}-a_{12}  \\  -a_9+a_{10}-a_{11}+a_{12}   & -a_5-a_6+a_7+a_8  &  a_1+a_2-a_3-a_4        \end{matrix}\right) $$
and 
$$ B= \left( \begin{matrix}   b_9+b_{10}-b_{11}-b_{12}  & -b_1+b_2-b_3+b_4 &  -b_5+b_6+b_7-b_8\\ b_1-b_2-b_3+b_4   & b_5+b_6-b_7-b_8 & b_9-b_{10}+b_{11}-b_{12} \\ -b_5+b_6-b_7+b_8 & b_{9}-b_{10}-b_{11}+b_{12} & -b_1-b_2+b_3+b_4  \end{matrix} \right)    .$$
So for $G=S_4$  the group determinant takes the form
$$ \mathscr{D}_{G}\left(a_1,\ldots ,a_{12},b_{1},\ldots ,b_{12}\right) =\ell_1 \ell_2 \;q_1^2 \;d_1^3d_2^3. $$

We can think of the determinant as the unormalized Lind  measure of the `polynomial' (really an element in the group ring $\mathbb Z [S_4]$)
\begin{align*} a_1 & +a_2x^2+a_3yx^2y+a_4x^2yx^2y+a_5xy+a_6x^3y +a_7xyx^2+a_8x^3yx^2 \\ & +a_9yx^3+a_{10}x^2yx+a_{11}x^2yx^3+a_{12}yx + b_1x+b_2x^3+b_3 xyx^2y + b_4yx^2yx \\ & + b_5x^3yx + b_6xyx^3+b_7xyx + b_8 x^3yx^3+ b_9y+b_{10}x^2yx^2+b_{11}yx^2+b_{12}x^2y, 
\end{align*}
where  monomials do not commute but we can reduce a polynomial in $\mathbb Z[x,y]$ to this form using the group relations $y^2=1$, $x^4=1$, $yxy=x^3yx^3$, $yx^3y=xyx$, $yx^2yx^2=x^2yx^2y$, $x^3yx^2y=yx^2yx,$ $x^2yx^2yx=xyx^2y$ etc.

\section{Proof of Theorem \ref{main}}
\begin{proof}
We first show that we can achieve the stated values. 

We begin with the values coprime to $3$ that are $1$ mod $4$.

Taking $a_1=1+k$ and the remaining values equal to $k$ gives 
$$\mathscr{D}_G=(1+24k)\cdot 1 \cdot (1-0)^2\cdot |I_3|^3\cdot |I_3|^3=1+24k. $$

Taking $a_2=a_5=a_9=1+k$, $b_3=b_5=1+k$ and the remaining values $k$ we get
$$ \mathscr{D}_G=(5+24k)\cdot 1 \cdot (0-1)^2\cdot \left| \begin{matrix} -1 & 0 & -2 \\ 0 & 0 & -1 \\ -2 & -1 & 2 \end{matrix}\right|^3  \cdot   \left|\begin{matrix} -1 & 2 & 0 \\ 2 & -2 & -1 \\ 0 & -1 & 0 \end{matrix}\right|^3=5+24k. $$

With $a_1=a_3=1+k$, $a_5=a_6=a_7=1+k$, $a_9=a_{10}=1+k$, and $b_1=b_3=1+k$, $b_5=b_6=1+k$, $b_{11}=b_{12}=1+k$ and the other values $k$ we obtain
$$\mathscr{D}_G=(13+24k)\cdot 1 \cdot (1-0)^2 \cdot \left| \begin{matrix}   -2 & 0 & -1 \\ -1 & 4 & 0 \\ 0 & -1 & 0 \end{matrix} \right|^3 \cdot    \left| \begin{matrix}   2 & 4 & -1 \\ -1 & 0 & 0 \\ 0 & -1 & 0 \end{matrix} \right|^3 =13+24k. $$

From $a_1=a_2=a_3=1+k$, $a_6=a_7=a_8=1+k$, $a_9=a_{10}=a_{11}=1+k$ and $b_1=b_3=b_4=1+k$, $b_5=b_7=b_8=1+k$, $b_9=b_{10}=1+k$, the others $k,$ we get
$$\mathscr{D}_G=(17+24k)\cdot 1 \cdot (0-1)^2 \cdot \left| \begin{matrix}   1 & 0 & 0 \\ 0 & 0 & 1 \\ -2 & 1 & 2 \end{matrix} \right|^3 \cdot    \left| \begin{matrix}   -3 & 2 & 2 \\ -2 & 2 & 1 \\ 0 & 1 & 0 \end{matrix} \right|^3 =17+24k. $$

Next we  get  the powers of three, $\pm 3^j$ with $j\geq 3$, that are $1$ mod $4$. 

With $a_1=a_3=1+k$, $b_3=1+k$ and the others $k$ we have
$$\mathscr{D}_G=(3+24k)\cdot 1 \cdot (4-1)^2 \cdot \left| \begin{matrix}   0 & -1 & 0 \\ -1 & 2 & 0 \\ 0 & 0 & 1 \end{matrix} \right|^3 \cdot    \left| \begin{matrix}   0 & 1 & 0 \\ 1 & 2 & 0 \\ 0 & 0 & -1 \end{matrix} \right|^3 =-3^3(1+8k), $$
and with $a_1=a_2=1+k$ and $a_5=1+k$ and the other values $k$
$$\mathscr{D}_G=(3+24k)\cdot 3 \cdot (3-0)^2 \cdot \left| \begin{matrix}   0 & 0 & -1 \\ 1 & 0 & 0 \\ 0 & -1 & 2 \end{matrix} \right|^6 =3^4(1+8k). $$

Finally we deal with the powers of two. 

Taking $a_1$, $a_5=1$ and the others zero  gives
$$\mathscr{D}_G=2\cdot 2 \cdot (1-0)^2 \cdot \left| \begin{matrix}   1 & 0 & -1 \\ 1 & 1 & 0 \\ 0& -1 & 1 \end{matrix} \right|^6 =2^8, $$
while $a_1=-1$, $a_5=a_6=1$ and $b_5=1$, $b_{10}=-1$ with the others zero  has
$$\mathscr{D}_G=1\cdot 1 \cdot (7-3)^2 \cdot \left| \begin{matrix}   -2 & 0 & -1 \\ 0 & 0 & 1 \\ -1 & -1 & -1 \end{matrix} \right|^3 \cdot    \left| \begin{matrix}   0 & 0 & 1 \\ 0 & -2 & -1 \\ 1 & -3 & -1 \end{matrix} \right|^3 =-2^{10}, $$
and $a_2=1$,$a_5=1$,$a_9=1$ and $b_{11}=1$ the rest zero
$$\mathscr{D}_G=4\cdot 2 \cdot (0-1)^2 \cdot \left| \begin{matrix}   -2 & 1 & -1 \\ 1 & -1 & 0 \\ -1 & -2 & 1 \end{matrix} \right|^3 \cdot    \left| \begin{matrix}   0 & 1 & -1 \\ 1 & -1 & -2 \\ -1 & 0 & 1 \end{matrix} \right|^3 =2^{12}, $$
and $a_1=a_2=1$, $a_5=-1$, $a_9=-1$, $b_1=1$ and the rest zero
$$\mathscr{D}_G=1\cdot -1  \cdot (9-1)^2 \cdot \left| \begin{matrix}  0 & -2 & 1 \\ 0 & 0 & 1 \\ 1 & 1 & 1 \end{matrix} \right|^3 \cdot    \left| \begin{matrix}   0 & 0 & 1 \\ -2 & 0 & 1 \\ 1 & 1 & 3 \end{matrix} \right|^3 =-2^{12}.$$
Taking $a_1=a_2=1+k$, $a_6=1+k$, $a_{10}=a_{11}=1+k$, $b_4=1+k$,  $b_6=1+k$, $b_{10}=1+k$
and the rest $k$
$$\mathscr{D}_G=(8+24k)\cdot 2  \cdot (1-0)^2 \cdot \left| \begin{matrix}  1 & 1 & 2 \\ 0 & 1 & 1 \\ 1 & -2 & 3 \end{matrix} \right|^3 \cdot    \left| \begin{matrix}   -1 & -1 & 0 \\ -2 & -1 & 3 \\ -1 & 0 & 1 \end{matrix} \right|^3 =2^{13}(1+3k),$$
and $a_2=a_3=a_4=1+k$, $a_5=1+k$,  $a_9=1+k$, $b_4=1+k$, $b_5=b_6=1+k$,  and the rest $k$
$$\mathscr{D}_G=(8+24k)\cdot 2  \cdot (4-3)^2 \cdot \left| \begin{matrix}  -1 & 2 & -1 \\ 2 & 1 & -1 \\ -1 & -1 & 0 \end{matrix} \right|^3 \cdot    \left| \begin{matrix}   -1 & 0 & -1 \\ 0 & -3 & -1 \\ -1 & -1 & -2 \end{matrix} \right|^3 =-2^{13}(1+3k),$$
giving  the remaining powers of two. Products of these achieve the stated values.

It remains to show that the determinants can only take the values claimed.  Plainly we have the congruences
\be \label{basic}   \ell_1\equiv \ell_2 \text{ mod } 2,\;\;\;  d_1\equiv d_2 \text{ mod } 2,\;\;\;   q_1\equiv \ell_1\ell_2 \text{ mod } 3, \ee
and, as can be checked on Maple,
\be \label{productmod4} d_1d_2 \equiv \ell_1\ell_2 q_1^2 \text{ mod } 4 \ee
and
\be \label{summod4}   d_1+d_2 \equiv (\ell_1+\ell_2 )q_1 \text{ mod } 4. \ee

Suppose that $3\mid \mathscr{D}_G=\ell_1\ell_2\; q_1^2\; (d_1d_2)^3.$ If $3\mid d_1d_2$ then  $3^3\mid (d_1d_2)^3$. If  $3\mid \ell_1\ell_2$  or $3\mid q_1$
then  $3$ divides both by \eqref{basic} and $3^3\mid \ell_1\ell_2 q_1^2$. Hence $3\nmid \mathscr{D}_G$
or $3^3\mid \mathscr{D}_G$.

If $\mathscr{D}_G$ is odd then by \eqref{productmod4} we have $\mathscr{D}_G\equiv (\ell_1\ell_2 q_1^2)^4\equiv 1$ mod $4$.

Suppose that $\mathscr{D}_G$ is even. From \eqref{productmod4} we have $2\mid d_1d_2$ and $2\mid \ell_1\ell_2q_1$.
Hence by \eqref{basic} we have $2\mid d_1$ and $2\mid d_2,$ and either $2\mid q_1$ or $2\mid \ell_1$ and $2\mid \ell_2$.
Hence $2^2\mid (\ell_1\ell_2)q_1^2$, $2^6\mid  (d_1d_2)^3$ and $2^8\mid \mathscr{D}_G$. It remains to rule out $2^9\parallel \mathscr{D}_G$ and $2^{11}\parallel \mathscr{D}_G$, and show that the $2^8\parallel \mathscr{D}_G$ or $2^{10}\parallel \mathscr{D}_G$
are of the stated forms.

As can be checked on Maple, expanding the determinant $d_1$ gives
\be \label{d1mod4}  d_1=\ell_1 \left( q_1 +2uv+ 2w\right)+4C(\vec{a},\vec{b}) \ee
where $C(\vec{a},\vec{b})$ is a homogeneous cubic integer  polynomial  in the $a_1,\ldots ,a_{12}$,
$b_{1},\ldots ,b_{12}$,
$$ u:=u_1+u_2+u_3, \;\;\; v:=v_1+v_2+v_3, $$
and 
$$ w:=u_1B_1+u_2B_2+u_3B_3+v_1A_1+v_2A_2+v_3A_3, $$
where
\begin{align*}
A_1:&= a_1+a_2+a_5+a_8+a_{9}+a_{10}, & B_1 & = b_1+b_2+b_7+b_8+b_{11}+b_{12}, \\
A_2 :&=a_1+a_3+a_6+a_8+a_{10}+a_{12}, &  B_2&=b_2+b_3+b_5+b_8+b_{10}+b_{11}, \\
A_3:&= a_1+a_4+a_7+a_8+a_{10}+a_{11}, & B_3& =b_1+b_3+b_5+b_7+b_{10}+b_{12}.\end{align*}
Similarly, replacing the the $b_i$ by $-b_i$
$$ d_2=\ell_2 \left( q_1 -2uv-2w\right)+4C(\vec{a},-\vec{b}) $$
and we get
\begin{align} 
d_1+d_2&  \equiv 2uq_1   +4uv^2 +4vw \text{ mod } 8,\label{mod8} \\
d_1-d_2 & \equiv 2v q_1 +4u^2 v  +4uw \text{ mod } 8. \label{mod8minus}
\end{align}

Suppose that $2^2\mid \ell_1$  or $2\mid \ell_1$ and $2\mid q_1$, then from \eqref{d1mod4} we have $d_1\equiv 0$ mod $4$ and 
and, since  $2\mid \ell_2,d_2$, we get $2^{12}\mid \mathscr{D}_G$.
Similarly if $2^2\mid \ell_2$ using $d_2$. 
Hence if $2^{12}\nmid \mathscr{D}_G$ we can assume that $2\nmid \ell_1\ell_2$ and $2\mid q_1$ or $2\parallel \ell_1,\ell_2$ and $2\nmid q_1$.  If $2\parallel \ell_1,\ell_2$ then $4\mid (\ell_1+\ell_2)$
and if $2\nmid \ell_1\ell_2$ then $2\mid q_1$  and $2\mid (\ell_1+\ell_2)$. In either case \eqref{summod4} gives 
$d_1\equiv d_2$ mod $4$, and  if $2^2\mid d_1$ or $d_2$ then $2^2$ divides both and $2^{14}\mid \mathscr{D}_G$. Hence we have $2\parallel d_1,d_2$ and so an even number of 2's divide $(d_1d_2)^2$, $q_1^2$ and $\ell_1\ell_2$, ruling out $2^9$ or $2^{11}\parallel \mathscr{D}_G$.

Suppose that $2^8$ or $2^{10}\parallel \mathscr{D}_G$. So $d_1=2\delta_1$, $d_2=2\delta_2$ with $\delta_1,\delta_2$ odd, and for $2^{10}\mid \mathscr{D}_G$ we must have $2\nmid \ell_1\ell_2$ and $2^2\parallel q_1$, and for  $2^8\parallel \mathscr{D}_G$ we have either $2\nmid \ell_1,\ell_2$ and $2\parallel q_1$, or $2\parallel \ell_1,\ell_2$ and $2\nmid q_1$. Suppose that $2\nmid \ell_1,\ell_2$. If $\ell_1\ell_2=u^2-v^2\equiv 1$  mod $4$  then $2\nmid u$, $2\mid v$ and from \eqref{mod8} 
we have $\delta_1+\delta_2 \equiv uq_1$ mod $4$. For $\mathscr{D}_G=2^{10}m$ we have $q_1=2^2\rho$ with $\rho$ odd, and $\delta_1\equiv -\delta_2$ mod $4$ and $m=\ell_1\ell_2 \rho^2 (\delta_1\delta_2)^3 \equiv 1\cdot 1\cdot (-1)^3 \equiv -1$ mod $4$.
For $\mathscr{D}_G=2^8m$ we have $q_1=2\rho$ with $\rho$ odd and $\delta_1+\delta_2 \equiv 2$ mod $4$ and 
$\delta_1\equiv \delta_2$ mod $4$ and $m=\ell_1\ell_2 \rho^2 (\delta_1\delta_2)^3 \equiv 1\cdot 1\cdot 1^3 \equiv 1$ mod $4$.
Similarly when $\ell_1\ell_2\equiv -1$ mod $4$ we have $2\mid u$, $2\nmid v$, and  use \eqref{mod8minus} to get $\delta_1\equiv \delta_2$ mod $4$ for $2^{10}$
and $\delta_1\equiv  -\delta_2$ mod $4$ for $2^8$, and $m\equiv -1$ and $1$ mod $4$ respectively as before.

This leaves the case $\mathscr{D}_G=2^8m$ where $\ell_1=2\xi_1$, $\ell_2=2\xi_2$ with $\xi_1$,$\xi_2$ odd. If $\xi_1\xi_2\equiv 1$ mod $4$ then $4\mid v=(\xi_1-\xi_2)$ and $2\parallel u=(\xi_1+\xi_2)$ and from \eqref{mod8minus} we get $\delta_1\equiv \delta_2$ mod $4$ and $m=\xi_1\xi_2 q_1^2 (\delta_1\delta_2)^3 \equiv 1$ mod $4$. Likewise if $\xi_1\xi_2\equiv -1$ mod $4$
then $4\mid u$ and $2\parallel v$ and from \eqref{mod8} we have $\delta_1\equiv -\delta_2$ and again $m\equiv 1$ mod $4$.
\end{proof}

\end{document}